\documentclass[11pt]{article}%
\usepackage{amsfonts}
\usepackage{amsmath}
\usepackage{amssymb}
\usepackage{graphicx}%
\newtheorem{theorem}{Theorem}

\newtheorem{corollary}[theorem]{Corollary}

\newtheorem{example}[theorem]{Example}

\newtheorem{lemma}[theorem]{Lemma}

\newtheorem{proposition}[theorem]{Proposition}
\newtheorem{remark}[theorem]{Remark}

\begin{document}

\title{On the topological type of anticonformal square roots of automorphisms of even
order of Riemann surfaces}
\author{Antonio F. Costa}
\maketitle

\begin{abstract}
Let $S$ be a (compact)\ Riemann surface of genus greater than one. Two
automorphism of $S$ are topologically equivalent if they are conjugated by a
homeomorphism. The topological classification of automorphisms is a classical
problem and its study was initiated by J. Nielsen who in the thirties
classified conformal ones. The case of anticonformal automorphisms is more
involved and was solved by K. Yocoyama in the 80s-90s. In order to decide
whether two anti-conformal automorphisms are equivalent, it is usually
necessary to take into account many invariants, some of which are difficult to
compute. In this work we present some situations where the topological
equivalence is mainly due to the genus of some quotient surfaces and the
algebraic structure of the automorphism group.

An anticonformal square root of a conformal automorphism $f$ is an
anticonformal automorphism $g$ such that $g^{2}=f$ . Let $g_{1}$ and $g_{2}$
be anticonformal square roots of the same conformal automorphism of order $m$,
where $m$ is an even integer. If genus of $S/\left\langle g_{1},g_{2}%
\right\rangle $ is even and genus of $S/\left\langle g_{i}\right\rangle $ is
$\neq2$ we prove that $\left\langle g_{1}\right\rangle $ and $\left\langle
g_{2}\right\rangle $ are topologically equivalent. If genus of $S/\left\langle
g_{1},g_{2}\right\rangle $ is odd and $\left\langle g_{1},g_{2}\right\rangle $
is abelian we obtain that $\left\langle g_{1}\right\rangle $ and $\left\langle
g_{2}\right\rangle $ are topologically equivalent. We give examples to justify
the condition genus of $S/\left\langle g_{i}\right\rangle $ $\neq2$ and
$\left\langle g_{1},g_{2}\right\rangle $ abelian in each case.

\end{abstract}

{\small MSC 2020: 14H37, 30F10, 57M60}

\section{Introduction}

Riemann surfaces with automorphisms (conformal and anticonformal) appear
naturally when studying algebraic curves with equations having a special form.
For example, hyperelliptic curves, which admit an automorphism of order 2 with
quotient space the sphere, or curves defined by real polynomial equations,
which admit an anticonformal automorphism of order two (conjugation
involution). A first step in the study of automorphisms of Riemann surfaces is
the determination of their topological equivalence class as maps between
surfaces. Two automorphisms $f_{1}$ and $f_{2}$ of a Riemann surface $S$ are
said to be topologically equivalent if there exists a homeomorphism
$h:S\rightarrow S$ such that $f_{1}=h\circ f_{2}\circ h^{-1}$, in a similar
way two groups of automorphisms $G_{1}$ and $G_{2}$ are topologically
equivalent if there is a homeomorphism $h$ such that $G_{1}=h\circ G_{2}\circ
h^{-1}$. The topological types of the automorphism groups of Riemann surfaces
produce the equisymmetric stratification of the moduli space \cite{Br}, and
the starting point of the study of this stratification is the topological
classification of cyclic groups of automorphisms.

\bigskip

In this paper we will restrict our attention to compact Riemann surfaces of
genus greater than one, in such case the possible automorphisms are of finite
order. If $f$ is a conformal automorphism of a Riemann surface $S$, by a
classical result of J. Nielsen [N] (see also [Y1]), the topological type of
$f$ is completely determined by the order of $f$ and the local invariants at
the branch points of the covering $\pi:S\rightarrow S/\left\langle
f\right\rangle $. Note that with such invariants the genus of $S/\left\langle
f\right\rangle $ is determined. The local invariants are as follows, if $f$
has order $m$ and $p\in S/\left\langle f\right\rangle $ is a branch point with
branch index $n$ then $f^{m/n}$ is an automorphism of $S$, $\pi^{-1}(p)$ is a
set of fixed points and $f^{m/n}$ acts in open neigbourhood discs around each
point of $\pi^{-1}(p)$ topologically as a rotation of angle $2\pi\frac{q}{n}$,
where $q$ and $n$ are coprimes. The local invariants for the branch point $p$
are the branch index $n$ and the rotation number $q$. Remark that if
$h:S\rightarrow S$ is an orientation reversing homeomorphism the rotation
numbers of $h\circ f\circ h^{-1}$ are the inverse of the ones of $f$.

\bigskip

The topological classification of anticonformal automorphisms is more involved
than for conformal ones (see \cite{Y2}, \cite{Y3} and \cite{C1}). The
classification of anticonformal involutions (order two automorphisms) is
simple. If $i_{1}$, $i_{2}$ are two anticonformal involutions of $S$, $i_{1}$,
$i_{2}$ are topologically equivalent if and only if $S/\left\langle
i_{1}\right\rangle $ and $S/\left\langle i_{2}\right\rangle $ are homeomorphic
surfaces, i.e. they have the same genus, are both orientable or not and their
boundaries have the same number of connected components.

\bigskip

For anticonformal automorphisms of order $2m$, $m>1$, we need more invariants.
If $g_{1}$ and $g_{2}$ are anticonformal automorphisms, to be topologically
equivalent, the first conditions are that $g_{1}$ and $g_{2}$ have the same
order $2m$ and that $S/\left\langle g_{1}\right\rangle $ and $S/\left\langle
g_{2}\right\rangle $ are homeomorphic, i.e. are both orientable or
non-orientable surfaces, with the same genus and the same number of connected
components of the boundary. Other invariants are the rotation numbers for
fixed points of the even powers of the automorphisms $g_{i}$, these rotation
numbers are defined up inversion if $S/\left\langle g_{1}\right\rangle $ is
not orientable. Furthermore there are also some rotation numbers for the
curves preimage of boundary components of $S/\left\langle g_{i}\right\rangle $
which also yield topological invariants. Finally there are some homological
invariants to be considered in case the surfaces $S/\left\langle
g_{i}\right\rangle $ are non-orientable (see 3.1.2).

\bigskip

It is known that in the case where $g$ has order $2m$ with $m$ odd and greater
than $1$, there are examples where there are two square roots $g_{1}$ and
$g_{2}$ of a conformal automorphism, i.e. $g_{i}^{2}=f$, which are not
topologically equivalent (see examples 2.2 and 2.3 in \cite{C2}). In these
examples either $S/\left\langle g_{1}\right\rangle $ and $S/\left\langle
g_{2}\right\rangle $ are not homeomorphic or the rotation numbers on curves
preimage of boundary components of $S/\left\langle g_{1}\right\rangle $ and
$S/\left\langle g_{2}\right\rangle $ are different for $g_{1}$ and $g_{2}$.
Note that the rotation numbers of fixed points of the even powers of $g_{i}$
(that are also powers of $f$) in case of square roots automatically match.

\bigskip

In this article we are going to study whether the topological type of a cyclic
group generated by an anticonformal automorphism $g$ of order $2m$, with $m$
even, is determined by $g^{2}$. If $g_{1}$ and $g_{2}$ are anticonformal
automorphisms with $g_{i}^{2}=f$, the condition that $f$ has even order
implies that $S/\left\langle g_{1}\right\rangle $ and $S/\left\langle
g_{2}\right\rangle $ are homeomorphic and have no boundary, then no rotation
numbers on boundary components to be considered (as we have remarked, the
rotation numbers on branched points are equivalent for $g_{1}$ and $g_{2}$).
This opens up the possibility that the topological type of $\left\langle
g_{i}\right\rangle $ is completely given by $f$. In this paper we present
situations when the anticonformal roots are topological equivalent and
examples where this does not happen, which contradicts a result of \cite{C2}
(in this paper a result of \cite{C1} is used which has been corrected recently
in \cite{C3}).

\bigskip

The following theorem summarises the main results of this article:

\begin{theorem}
\label{Thm}Let $g_{1}$ and $g_{2}$ be two anticonformal automorphisms of $S$
such that $g_{1}^{2}=g_{2}^{2}=f$ and $f$ has even order $m$. Assume that the
order of $\left\langle g_{1},g_{2}\right\rangle $ is $n\times2m$.

a. If $n$ is even, genus of $S/\left\langle g_{1},g_{2}\right\rangle $ is even
and genus of $S/\left\langle g_{i}\right\rangle $ is $\neq2$ then $g_{1}$ and
$g_{2}$ are topologically equivalent.

b. If $n$ is even, genus of $S/\left\langle g_{1},g_{2}\right\rangle $ is odd
and $\left\langle g_{1},g_{2}\right\rangle $ is abelian then $g_{1}$ and
$g_{2}$ are topologically equivalent.

c. If $n$ is odd then $\left\langle g_{1}\right\rangle $ and $\left\langle
g_{2}\right\rangle $ are topologically equivalent.
\end{theorem}

\bigskip

Furthermore, we give examples to justify the condition genus of
$S/\left\langle g_{i}\right\rangle $ $\neq2$, in point a, and the condition
$\left\langle g_{1},g_{2}\right\rangle $ abelian, in point b of Theorem
\ref{Thm} (Examples \ref{Ex1} and \ref{Ex2}).

\bigskip

Note that the conditions in the theorem may be simpler to verify in order to
conclude the topological equivalence of $\left\langle g_{1}\right\rangle $ and
$\left\langle g_{2}\right\rangle $ than the calculation of the topological
invariants in \cite{Y2}, \cite{Y3} or \cite{C1}, \cite{C3}. As a conclusion it
can be said that two anticonformal square roots of an automorphism of even
order rarely have different topological types.

\bigskip

Finally, the results of this article can be used in the study of the branch
loci of the moduli space of Riemann surfaces. Two groups of automorphisms are
topologically equivalent if they are conjugate by an homeomorfism of the
surface. The topological types of automorphisms groups provide a
stratification of the moduli space given by the orbifold structure originated
from the covering $\mathbb{T}_{g}\rightarrow\mathcal{M}_{g}$. The strata
defined by cyclic groups are particularly important, since their closure cover
all the moduli space $\mathcal{M}_{g}$. Note that if two cyclic automorphism
groups are topologically conjugate then there are generators of such groups
that must be topologically equivalent. The results of this paper indicate that
stratification by topological types of groups of anticonformal automorphisms
can distinguish points belonging to the closure of a single stratum of
conformal automorphisms of a given topological type of even order.

The surfaces with anticonformal automorphisms appear in $\mathcal{M}_{g}$ as
those points which are fixed by the automorphism conjugation $-:\mathcal{M}%
_{g}\rightarrow\mathcal{M}_{g}$, and constitute the set of real and pseudoreal
surfaces. The pseudoreal surfaces are surfaces admitting anticonformal
automorphisms but without anticonformal involutions, see \cite{BCC}, then
having anticonformal automorphisms of order $2m$, with $m$ even. The
topological types of these automorphisms provide equisymmetric strata whose
closure covers the set of pseudoreal Riemann surfaces.

\bigskip

\section{Preliminaries on automorphisms of Riemann surfaces}

\subsection{Non-Euclidean crystallographic groups}

A non-Euclidean crystallographic (NEC) group $\Lambda$ is a group of conformal
and anticonformal automorphisms of the complex upper plane with a compact
2-dimensional orbifold as orbit space $\mathbb{H}^{2}/\Lambda$. The group
$\Lambda$ is isomorphic to the orbifold fundamental group $\pi_{1}%
O(\mathbb{H}^{2}/\Lambda)$.

The signature of $\Lambda$ is a symbol as follows:
\[
(g;\pm;[m_{1},...,m_{r}],\{(n_{11},...,n_{1s_{1}}),...,(n_{k1},...,n_{ks_{k}%
})\})
\]
where $(g,\pm,k)$ is the genus, the orientability and the number of connected
boundary components respectively of $\mathbb{H}^{2}/\Lambda$ as a topological
surface, then $(g,\pm,k)$ determines the topology of $\mathbb{H}^{2}/\Lambda$.
The orbifold $\mathbb{H}^{2}/\Lambda$ has $r$ conic points with isotropy
cyclic groups of orders $m_{1},...,m_{r}$, and each boundary connected
component has $s_{j}$ corner points, the isotropy groups of the corner points
are dihedral groups $D_{n_{ij}},$ where $1\leq i\leq k$ and $1\leq j\leq
s_{j}$. The signature describes completely the algebraic structure of the
group $\Lambda$ and the orbifold structure of $\mathbb{H}^{2}/\Lambda$. The
permutations in the set of periods (numbers) $[m_{1},...,m_{r}]$, the cyclic
permutations in the numbers in the brackets $(n_{11},...,n_{1s_{1}})$ and
permutations of the order of the brackets, produce equivalent signatures (see
\cite{BEGG} or \cite{W}).

The group $\Lambda$ has canonical presentations as follows. If the sign in the
signature is $+$, then:%
\[
\left\langle
\begin{array}
[c]{c}%
a_{1},b_{1},...,a_{g},b_{g},x_{1},...,x_{r},e_{1},...,e_{k},c_{ij},0\leq i\leq
s_{j},1\leq j\leq k:\\%
{\textstyle\prod\nolimits_{l=1}^{g}}
[a_{l},b_{l}]%
{\textstyle\prod\nolimits_{l=1}^{r}}
x_{l}%
{\textstyle\prod\nolimits_{l=1}^{k}}
e_{l}=1;x_{l}^{m_{l}}=1;c_{ij}^{2}=1;e_{l}c_{ls_{l}}e_{l}^{-1}=c_{l0}%
\end{array}
\right\rangle
\]
if the sign in the signature is $-$ and $g$ is odd:%
\[
\left\langle
\begin{array}
[c]{c}%
d,a_{1},b_{1},...,a_{\frac{g-1}{2}},b_{\frac{g-1}{2}},x_{1},...,x_{r}%
,e_{1},...,e_{k},c_{ij},0\leq i\leq s_{j},1\leq j\leq k:\\
d^{2}%
{\textstyle\prod\nolimits_{l=1}^{\frac{g-1}{2}}}
[a_{l},b_{l}]%
{\textstyle\prod\nolimits_{l=1}^{r}}
x_{l}%
{\textstyle\prod\nolimits_{l=1}^{k}}
e_{l}=1;x_{l}^{m_{l}}=1;c_{ij}^{2}=1;e_{l}c_{ls_{l}}e_{l}^{-1}=c_{l0}%
\end{array}
\right\rangle
\]
and if the sign in the signature is $-$ and $g$ is even:%
\[
\left\langle
\begin{array}
[c]{c}%
d_{1},d_{2},a_{1},b_{1},...,a_{\frac{g-1}{2}},b_{\frac{g-1}{2}},x_{1}%
,...,x_{r},e_{1},...,e_{k},c_{ij},0\leq i\leq s_{j},1\leq j\leq k:\\
d_{1}d_{2}%
{\textstyle\prod\nolimits_{l=1}^{\frac{g-2}{2}}}
[a_{l},b_{l}]%
{\textstyle\prod\nolimits_{l=1}^{r}}
x_{l}%
{\textstyle\prod\nolimits_{l=1}^{k}}
e_{l}=1;x_{l}^{m_{l}}=1;c_{ij}^{2}=1;e_{l}c_{ls_{l}}e_{l}^{-1}=c_{l0}%
\end{array}
\right\rangle
\]

The generators $a_{i},b_{i}$ are hyperbolic transformations and $d,d_{1}%
,d_{2}$ are glide reflections. The $c_{ij}$ are reflections, the $x_{i}$ are
elliptic transformations and the $e_{i}$ are, in the general case, hyperbolic
transformations. The generators of a canonical presentation of an NEC\ group
are called canonical generators.

If the group $\Lambda$ does not contain anticonformal transformations then
$\Lambda$ is a Fuchsian group and its signature has the form $(g;+;[m_{1}%
,...,m_{r}])$, i. e. Fuchsian groups can be considered as a special type of NEC\ groups.

In the study of the topology of automorphisms we need to consider the
abelianization $\Lambda/[\Lambda,\Lambda]$ of a NEC\ group $\Lambda$. The
group $\Lambda/[\Lambda,\Lambda]$ will be denoted by $H_{1}O(\mathbb{H}%
^{2}/\Lambda)$, the first orbifold homology group of $\mathbb{H}^{2}/\Lambda$.
If $\alpha:\Lambda\rightarrow\Lambda/[\Lambda,\Lambda]=H_{1}O(\mathbb{H}%
^{2}/\Lambda)$ is the abelianization morphism we shall denote by capital
letters the images by $\alpha$ of the generators of a canonical presentation
of $\Lambda$:
\[
\alpha(a_{i})=A_{i},\alpha(b_{i})=B_{i},\alpha(d_{i})=D_{i},\ \alpha
(x_{i})=X_{i},\alpha(e_{i})=E_{i},\alpha(c_{i})=C_{i}.
\]

\subsection{Fuchsian and NEC\ groups and automorphisms of Riemann surfaces}

A Fuchsian \textit{surface} group is a Fuchsian group with signature
$(g,+,[-];\{-\})$. Let $S$ be a Riemann surface of genus $g\geq2$, by the
Uniformization Theorem of Riemann surfaces, there is a Fuchsian surface group
$\Gamma$ such that $S=\mathbb{H}^{2}/\Gamma$. Let $\mathrm{Aut}^{\pm}(S)$ be
the group of conformal and anticonformal automorphisms $S$, if $G$ is a
subgroup of $\mathrm{Aut}^{\pm}(S)$ then there is an NEC group $\Delta$ such
that $\Gamma\vartriangleleft\Delta$, $G\cong\Delta/\Gamma$ and $S/G=\mathbb{H}%
^{2}/\Delta$. If $G$ is contained in $\mathrm{Aut}^{+}(S)$ ($=\mathrm{Aut}%
(S)$) (there is no anticonformal automorphisms in $G$) the group $\Delta$ is a
Fuchsian group.

If $\left\langle f\right\rangle $ is a cyclic group generated by a conformal
automorphism of degree $m$ the group $\Delta$ containing $\Gamma$, such that
$S/\left\langle f\right\rangle =\mathbb{H}^{2}/\Delta$, is a Fuchsian group
with signature $(q;+;[m_{1},...,m_{r}])$, where the $m_{i}$ are divisors of
$m$. Let us call $p:S\rightarrow S/\left\langle f\right\rangle $ and
$\pi:\mathbb{H}^{2}\rightarrow\mathbb{H}^{2}/\Gamma$.

There is an epimorphism
\[
\omega_{f}:\Delta\rightarrow\Delta/\Gamma\cong\left\langle f\right\rangle
\cong\mathbb{Z}_{m}%
\]
given by the action of $f$. The definition of $\omega_{f}$ is as follows.
Consider a point $o$ in $S/\left\langle f\right\rangle $, and now choose a
point $\widetilde{o}_{0}\in p^{-1}(o)$ and denote by $\widetilde{o}_{i}%
=f^{i}(\widetilde{o}_{0})$ each one of the lifts of $o$ to $S$. Let $y$ be an
element of $\Delta$ is a transformation of $\mathbb{H}^{2}$. Let
$\widetilde{o}$ be a lift of $\widetilde{o}_{0}$ to $\mathbb{H}^{2}$. We
define $\omega_{f}(y)=\overline{i}\in\mathbb{Z}_{m}$ if $\pi(y(\widetilde{o}%
))=\widetilde{o}_{i}$. The epimorphism $\omega_{f}$ is the monodromy of the
action of $f$. Let $\Delta/[\Delta,\Delta]\cong H_{1}O(\mathbb{H}^{2}%
/\Delta)=H_{1}O(S/\left\langle f\right\rangle )$ be the orbifold first
homology group for $S/\left\langle f\right\rangle $. Since $\mathbb{Z}_{m}$ is
abelian the monodromy $\omega_{f}$ is determined by $\Omega_{f}:H_{1}%
O(S/\left\langle f\right\rangle )\rightarrow\mathbb{Z}_{m}$, where
\[
\omega_{f}:\pi_{1}O(S/\left\langle f\right\rangle )\cong\Delta\overset{\alpha
}{\rightarrow}H_{1}O(S/\left\langle f\right\rangle )\cong\Delta/[\Delta
,\Delta]\overset{\Omega_{f}}{\rightarrow}\mathbb{Z}_{m}.
\]

In a similar way we study anticonformal automorphisms. Assume $\left\langle
g\right\rangle $ a cyclic group of order $2m$ such that a generator $g$ is an
anticonformal automorphism. Since the dihedral isotropy groups cannot be
subgroups of $\left\langle g\right\rangle $ then the signature of $\Delta$ is
as follows:%
\[
(h;\pm;[m_{1},...,m_{r}],\{(-),\overset{k}{...},(-)\})
\]
where the $m_{i}$ are divisors of $m$. As in the case of conformal
automorphisms we define the epimorphisms $\omega_{g}:\Delta\rightarrow
\mathbb{Z}_{2m}$ and $\Omega_{g}:H_{1}O(S/\left\langle g\right\rangle
)\rightarrow\mathbb{Z}_{2m}$.

Moreover, if $m$ is an even integer the signature of $\Delta$ have the form
$(h;-;[m_{1},...,m_{r}])$, note that in this case the unique element of order
two in $\left\langle g\right\rangle $ is $g^{m}$ that is a conformal
automorphism, as consequence there is no reflections in $\Delta$.

\section{Topological preliminaries}

Two automorphisms \thinspace$g_{1}$ and $g_{2}$ of a Riemann surface $S$ are
topological equivalent (or have the same topological type)\ if and only if
there is a homeomorphism $h:S\rightarrow S$ such that $g_{1}=h\circ g_{2}\circ
h^{-1}$.

Let $g_{i}$, $i=1,2$, be two anticonformal automorphisms of order $2m$ of the
Riemann surface $S$, $\Gamma$ be a Fuchsian surface group uniformizing $S$ and
$\Delta_{i}$, $i=1,2$ be NEC\ crystallographic groups such that $\Gamma
\vartriangleleft\Delta_{i}$ and $S/\left\langle g_{i}\right\rangle
\cong\mathbb{H}^{2}/\Delta_{i}$. If $g_{1}$ is topologically equivalent to
$g_{2}$ then the signature of $\Delta_{1}$ is equivalent to the signature of
$\Delta_{2}$. By the Dehn-Nielsen Theorem and branched covering theory (see
\cite{S}) $g_{1}$ and $g_{2}$ are topologically equivalent if and only if
there is an isomorphism $\iota:\Delta_{1}\rightarrow\Delta_{2}$ such that
$\omega_{g_{1}}=\omega_{g_{2}}\circ\iota$ where $\omega_{g_{i}}:\Delta
_{i}\rightarrow\mathbb{Z}_{2m}$ are the monodromies of $g_{1}$ and $g_{2}$.
Equivalently $g_{1}$ and $g_{2}$ are topologically equivalent if and only if
there is an isomorphism $\iota:\Delta_{1}\rightarrow\Delta_{2}$ such that
$\Omega_{g_{1}}=\Omega_{g_{2}}\circ\iota_{\ast}$ where $\iota_{\ast}%
:H_{1}O(S/\left\langle g_{1}\right\rangle )\rightarrow H_{1}O(S/\left\langle
g_{2}\right\rangle )$ is the isomorphism induced on homology by $\iota$.

\bigskip

\subsection{Classification theorems}

Let $g_{1}$ and $g_{2}$ be two anticonformal automorphism of order $2m$ and
$\Omega_{g_{1}}:H_{1}O(S/\left\langle g_{1}\right\rangle )\rightarrow
\mathbb{Z}_{2m}$, $\Omega_{g_{2}}:H_{1}O(S/\left\langle g_{2}\right\rangle
)\rightarrow\mathbb{Z}_{2m}$, be the corresponding monodromies.

Assume that the quotient orbifolds $S/\left\langle g_{1}\right\rangle $ and
$S/\left\langle g_{2}\right\rangle $ have the same signature%
\[
(h;\pm;[m_{1},...,m_{r}];\{(-),\overset{k}{...},(-)\}).
\]

\subsubsection{Case 1. $S/\left\langle g_{t}\right\rangle $ orientable (then
$m$ is odd).}

Assume that $\{A_{i}^{(t)},B_{i}^{(t)},X_{j}^{(t)},E_{l}^{(t)}\}$ are the
generators of $H_{1}O(S/\left\langle g_{t}\right\rangle )$ given by a
canonical presentation of $\pi_{1}O(S/\left\langle g_{t}\right\rangle )$.

The automorphisms $g_{1}$ and $g_{2}$ are topologically equivalent if and only
if
\begin{align*}
&  (\{\Omega_{g_{1}}(X_{1}),...,\omega_{g}(X_{r})\},\{\Omega_{g_{1}}%
(E_{1}),...,\omega_{g}(E_{k})\})\\
&  =(\{\varepsilon\omega_{g}(X_{1}),...,\varepsilon\omega_{g}(X_{r}%
)\},\{\varepsilon\omega_{g}(E_{1}),...,\varepsilon\omega_{g}(E_{k})\})
\end{align*}
where $\varepsilon=\pm1$ (see \cite{Y3}, \cite{C1}).

\subsubsection{Case 2. $S/\left\langle g_{t}\right\rangle $ non-orientable.
\label{Classification}}

The automorphisms $g_{1}$ and $g_{2}$ are topologically equivalent if and only if:

There are generator systems $\{D_{q}^{(t)},A_{i}^{(t)},B_{i}^{(t)},X_{j}%
^{(t)},E_{l}^{(t)}\}_{t=1,2}$ of $H_{1}O(S/\left\langle g_{t}\right\rangle )$
given by canonical presentations of $\pi_{1}O(S/\left\langle g_{t}%
\right\rangle )$ such that the following three properties are verified:

1.
\[
(\{\Omega_{g_{1}}(X_{j}^{(1)})\},\{\Omega_{g_{1}}(E_{l}^{(1)})\})=(\{\Omega
_{g_{1}}(X_{j}^{(2)})\},\{\Omega_{g_{1}}(E_{l}^{(2)})\})
\]

2. If there is no order two elements in $\{\Omega_{g_{1}}(X_{j}^{(1)}%
)\},\{\Omega_{g_{1}}(E_{l}^{(1)})\}$:
\[
\Omega_{g_{1}}(D^{(1)})=\Omega_{g_{1}}(D^{(2)})
\]
for $h$ odd and%
\[
\Omega_{g_{1}}(D_{1}+D_{2})=\Omega_{g_{1}}(D_{1}+D_{2})
\]
for $h$ even.

3. If $h=2$:
\[
\Omega_{g_{1}}(D_{1})\equiv\pm\Omega_{g_{2}}(D_{1})\operatorname{mod}z
\]
where $z=\gcd\{\Omega_{g_{1}}(X_{j}^{(1)}),\Omega_{g_{1}}(E_{l}^{(1)}%
),\Omega_{g_{1}}(D_{1}+D_{2}),2m\}$.

\bigskip

\begin{remark}
\label{Remark h}
\end{remark}

Assume that there is an anticonformal automorphism $g$ and a generator system
$\{D_{q},A_{i},B_{i},X_{j},E_{l}\}$ of $H_{1}O(S/\left\langle g_{t}%
\right\rangle )$ given by a presentation of $\pi_{1}O(S/\left\langle
g\right\rangle )$. There is an automorphism $h_{D_{t},X_{j}}$ of $\pi
_{1}O(S/\left\langle g_{t}\right\rangle )$ such that
\[
\Omega_{g}((h_{D_{t},X_{j}})_{\ast}(D_{t}))=\Omega_{g}(D_{t})+\Omega_{g}%
(X_{j}),(\Omega_{g}((h_{D_{t},X_{j}})_{\ast}(X_{j}))=-\Omega_{g}(X_{j})
\]
and the rest of generators of $\{D_{q},A_{i},B_{i},X_{j},E_{l}\}$ are
invariant by $(h_{D_{t},X_{j}})_{\ast}$.

Similarly there is an automorphism $h_{D_{t},E_{j}}$ with:%
\[
\Omega_{g}((h_{D_{t},E_{j}})_{\ast}(D_{t}))=\Omega_{g}(D_{t})+\Omega_{g}%
(E_{j}),(\Omega_{g}((h_{D_{t},E_{j}})_{\ast}(E_{j}))=-\Omega_{g}(E_{j})
\]

\bigskip

The following Lemma will be useful in Section 4.

\begin{lemma}
\label{Lemma1}Let $g$ be an anticonformal automorphism of order $2m$ with $m$
even of a Riemann surface $S$. Let $\{A_{i},B_{i},D_{q},X_{j}\}$ be a
generator system of $H_{1}O(S/\left\langle g\right\rangle )$ given by a
canonical presentation of $\pi_{1}O(S/\left\langle g\right\rangle )$. Assume
that all the orders of the elements $X_{j}$ are different from $2$. Let $N$ be
a non-trivial homology class of $H_{1}O(S/\left\langle g\right\rangle )$ such
that $2N+%
{\textstyle\sum\nolimits_{j}}
\varepsilon_{j}X_{j}=0$, where $\varepsilon_{j}\in\{-1,1\}$.

Then there is a canonical presentation of $\pi_{1}O(S/\left\langle
g\right\rangle )$ producing
\[
\{A_{i}^{\prime},B_{i}^{\prime},D_{q}^{\prime},\varepsilon_{j}X_{j}\}
\]
as generator system of $H_{1}O(S/\left\langle g\right\rangle )$, such that $N$
is homologically equal to either $D_{1}^{\prime}$ (if the genus of
$S/\left\langle g\right\rangle $ is odd) or $D_{1}^{\prime}+D_{2}^{\prime}$
(if the genus of $S/\left\langle g\right\rangle $ is even).
\end{lemma}

\bigskip

\textbf{Proof.}

Since there are not $X_{j}$ of order $2$ and $N$ is non-trivial, the element
$N$ is completely determined by $2N+%
{\textstyle\sum\nolimits_{j}}
\varepsilon_{j}X_{j}=0$. Let $h$ be the automorphism of $\pi_{1}%
O(S/\left\langle g\right\rangle )$ defined by $h=%
{\textstyle\prod\nolimits_{\varepsilon_{j}=-1}}
h_{D_{1},X_{j}}$, where the $h_{D_{1},X_{j}}$ are the automorphisms in Remark
\ref{Remark h}. The image by $h$ of the canonical presentation giving
$\{A_{i},B_{i},D_{q},X_{j}\}$ is a presentation $\{A_{i}^{\prime}%
,B_{i}^{\prime},D_{q}^{\prime},\varepsilon_{j}X_{j}\}$ such that \ $%
{\textstyle\sum\nolimits_{q}}
2D_{q}^{\prime}+%
{\textstyle\sum\nolimits_{j}}
\varepsilon_{j}X_{j}=0$. then $N=%
{\textstyle\sum\nolimits_{q}}
D_{q}^{\prime}$. \ \ $\square$

\bigskip

\textbf{Note}.

Let $H$ be an homology class of $H_{1}O(S/\left\langle g\right\rangle )$ such
that the image of $H$ by $H_{1}O(S/\left\langle g\right\rangle )\rightarrow
H_{1}O(S/\left\langle g\right\rangle )/\left\langle X_{j}\right\rangle \cong
H_{1}(S/\left\langle g\right\rangle )$ is of order $2$ then, by \cite{MIII}
$H$ can be represented by a simple closed curve in $S/\left\langle
g\right\rangle $.

\bigskip

\subsection{Topological invariants}

Let us give a topological interpretation of conditions 2 and 3 of the above
classification in subsection \ref{Classification}. The condition 2 is related
with the existence of simple closed simple curves on non-orientable
2-orbifolds such that cutting along such curves we obtain orientable
orbifolds with one or two boundary components. More precisely, the existence
of curves $\gamma_{1}$ and $\gamma_{2}$ in $S/\left\langle g_{1}\right\rangle
$ and $S/\left\langle g_{2}\right\rangle $ respectively in such a way that
$S-p_{1}^{-1}(\gamma_{1})\rightarrow$ $S/\left\langle g_{1}\right\rangle
-\gamma_{1}$ is topologically equivalent to $S-p_{2}^{-1}(\gamma
_{2})\rightarrow$ $S/\left\langle g_{2}\right\rangle -\gamma_{2}$, as branched
coverings of orientable surfaces with boundary. With this condition the
homology of $[\gamma_{1}]$ and $[\gamma_{2}]$ are determined and the action of
$g_{1}$ on $p_{1}^{-1}(\gamma_{1})$ must be equivalent to the action of
$g_{2}$ on $p_{2}^{-1}(\gamma_{2})$. This is the topological meaning of
condition 2.

For condition three there is a fixed point free action of $g_{i}^{2m/z}$ on
$S/\left\langle g_{i}^{z}\right\rangle $ with quotient $S/\left\langle
g_{t}\right\rangle $, $t=1,2$. The quotient surfaces $S/\left\langle g_{t}%
^{z}\right\rangle $ are homeomorphic to the Klein bottle and in order to have
the equivalence of the free actions of $g_{1}^{2m/z}$ and $g_{2}^{2m/z}$ it
appears condition 3.

In the next example we show a case where the genus of $S/\left\langle
g\right\rangle $ is two and $g$ is not topologically equivalent to $g^{m+1}$.
Note that $g^{2}=(g^{m+1})^{2}$, so this is a simple case where, for $m$ even,
the topological type of $g$ is not determined by $g^{2}$.

\bigskip

Consider $G=C_{8}\times C_{3}=\left\langle u:u^{8}=1\right\rangle
\times\left\langle t:t^{3}=1\right\rangle $ and an NEC\ group $\Delta$ of
genus $2$ with signature $(2,-,[3,3])$. A canonical presentation is:%
\[
\left\langle d_{1},d_{2},x_{1},x_{2}:d_{1}^{2}d_{2}^{2}x_{1}x_{2}=1;x_{i}%
^{3}=1,i=1,2\right\rangle
\]

We construct the following monodromy:%
\begin{align*}
\omega &  :\Delta\rightarrow G\\
d_{1}  &  \longmapsto(u,1),d_{2}\longmapsto(u^{7},1),x_{1}\longmapsto
(1,t),x_{2}\longmapsto(1,t^{2})
\end{align*}

We consider $S=\mathbb{H}^{2}/\ker\omega$. The group $G$ acts on $S$ and the
elements $g_{1}=(u,t)$ and $g_{2}=(u^{5},t)$ have the same square in $G$. We
have $S/\left\langle g_{1}\right\rangle =\mathbb{H}^{2}/\omega^{-1}%
(\left\langle (u,t)\right\rangle )=\mathbb{H}^{2}/\omega^{-1}(\left\langle
(u^{5},t)\right\rangle )=S/\left\langle g_{2}\right\rangle $ and the genus of
$S/\left\langle g_{i}\right\rangle $ is two.

Let $\{D_{1},D_{2},X_{1},X_{2}\}$ be a generator system of $H_{1}%
O(S/\left\langle g_{t}\right\rangle )$ obtained from the above canonical
presentation of $\Delta$, i. e. $\alpha(d_{i})=D_{i}$, $\alpha(x_{j})=X_{j}$.

To apply the classification we need to compute some values of $\Omega_{g_{1}}$
and $\Omega_{g_{2}}$ on some generators:
\begin{gather*}
\omega(x_{1})=(1,t)=g_{1}^{8}=g_{2}^{8}\text{, then}\\
\Omega_{g_{1}}(X_{1})=\overline{8}=\Omega_{g_{2}}(X_{1})\\
\omega(x_{2})=(1,t^{2})=g_{1}^{16}=g_{2}^{16}\text{, then}\\
\Omega_{g_{1}}(X_{2})=\overline{16}=\Omega_{g_{2}}(X_{2})\\
\omega(d_{1}d_{2})=1_{G}=\omega(d_{1}d_{2})\\
\Omega_{g_{1}}(D_{1}+D_{2})=\overline{0}=\Omega_{g_{2}}(D_{1}+D_{2})
\end{gather*}

Finally $\Omega_{g_{1}}(D_{1})=\overline{9}$ and $\Omega_{g_{2}}%
(D_{1})=\overline{21}$ and the number $z\ $is
\[
\gcd\{\Omega_{g_{1}}(X_{j}),\Omega_{g_{1}}(D_{1}+D_{2}),24\}=8.
\]
Since
\[
\Omega_{g_{1}}(D_{1})\operatorname{mod}8\equiv1\operatorname{mod}8\neq
\Omega_{g_{2}}(D_{1})\operatorname{mod}8\equiv5\operatorname{mod}8,
\]
then $g_{1}$ and $g_{2}$ are not topologically equivalent. $\square$

\bigskip

\section{Anticonformal square roots of a conformal automorphism of order
even.}

We begin by setting some notations that will be used in this section. Let $S$
be a Riemann surface and $f$ a conformal automorphism of $S$ of even order
$m$. Let $g_{1}$ and $g_{2}$ be two anticonformal automorphisms such that
$g_{i}^{2}=f$, $i=1,2$ ($g_{i}$ are square roots of $f$).

The group $\left\langle f\right\rangle $ is a normal subgroup of $\left\langle
g_{1},g_{2}\right\rangle $ and $\left\langle g_{1},g_{2}\right\rangle
/\left\langle f\right\rangle $ is a dihedral group. The involutions generating
$\left\langle g_{1},g_{2}\right\rangle /\left\langle f\right\rangle $ do not
have fixed curves, then the projection $S/\left\langle f\right\rangle
\rightarrow S/\left\langle g_{1},g_{2}\right\rangle $ is a branched covering
with deck transformations group a dihedral group $D_{n}$, i. e,. $\left\langle
g_{1},g_{2}\right\rangle /\left\langle g_{i}^{2}\right\rangle =\left\langle
g_{1},g_{2}\right\rangle /\left\langle f\right\rangle \cong D_{n}$, where $n$
is an integer. We denote
\[
\theta:\left\langle g_{1},g_{2}\right\rangle =G\rightarrow\left\langle
g_{1},g_{2}\right\rangle /\left\langle f\right\rangle =D_{n}=\left\langle
a,b:a^{2}=b^{2}=(ab)^{2}=1\right\rangle ,
\]
defined by $\theta(g_{1})=a$; $\theta(g_{2})=b$.

\bigskip

\begin{proposition}
Let $S$ be a Riemann surface and $f$ be a conformal automorphism of $S$ of
order $m$. Let $g_{1}$ and $g_{2}$ be two anticonformal automorphisms such
that $g_{i}^{2}=f$. If the order of $\left\langle g_{1},g_{2}\right\rangle $
is $2m\times n$, where $n$ is odd, then $\left\langle g_{1}\right\rangle $ and
$\left\langle g_{2}\right\rangle $ are conformally conjugate.
\end{proposition}

\textbf{Proof. }

If $\theta:\left\langle g_{1},g_{2}\right\rangle \rightarrow\left\langle
g_{1},g_{2}\right\rangle /\left\langle f\right\rangle =D_{n}$, since $n$ is
odd, $\theta(g_{1})$ and $\theta(g_{2})$ are conjugate in $D_{n}$. Then
\[
\theta(g_{1})=\theta((g_{1}g_{2})^{\alpha})\theta(g_{2})\theta((g_{1}%
g_{2})^{-\alpha})=\theta((g_{1}g_{2})^{\alpha}g_{2}(g_{1}g_{2})^{-\alpha}),
\]
then $g_{1}=(g_{1}g_{2})^{\alpha}g_{2}(g_{1}g_{2})^{-\alpha}g_{2}^{2w}%
=(g_{1}g_{2})^{\alpha}g_{2}^{2w+1}(g_{1}g_{2})^{-\alpha}$. Hence $\left\langle
g_{1}\right\rangle $ and $\left\langle g_{2}\right\rangle $ are conformally
conjugate. $\ \ \square$

\bigskip

Now we shall consider the case $n$ even.

By the Uniformization Theorem we have $S\cong\mathbb{H}^{2}/\Gamma$ with
$\Gamma$ a Fuchsian surface group and let $\Delta$ be a NEC\ group such that
$\Gamma\subset\Delta$ and $\mathbb{H}^{2}/\Delta\cong S/\left\langle
g_{1},g_{2}\right\rangle $. Since $m$ is even, there are no anticonformal
involutions in the group $\left\langle g_{i}\right\rangle $ and then
$S/\left\langle g_{i}\right\rangle $, $i=1,2$, are non-orientable surfaces
without boundary. The anticonformal involutions generating the deck
transformation group of $S/\left\langle f\right\rangle \rightarrow
S/\left\langle g_{1},g_{2}\right\rangle $ have no fixed points. Hence there
are no reflections in $\Delta$ and the NEC\ group $\Delta$ has signature
$(h;-;[m_{1},...,m_{r}];\{-\})$ where the integer $h$ is the genus of
$S/\left\langle g_{1},g_{2}\right\rangle $.

\bigskip

Let $\omega:\Delta\rightarrow G=\left\langle g_{1},g_{2}\right\rangle $ be the
monodromy of $S\rightarrow S/\left\langle g_{1},g_{2}\right\rangle $ then
$\ker\omega=\Gamma$. We call $\omega^{-1}\circ\theta^{-1}(a)=\Delta_{1}$ and
$\omega^{-1}\circ\theta^{-1}(b)=\Delta_{2}$. We have that $S/\left\langle
g_{1}\right\rangle =\mathbb{H}^{2}/\Delta_{1}$ and $S/\left\langle
g_{2}\right\rangle =\mathbb{H}^{2}/\Delta_{2}$. Note that the monodromy of
$S\rightarrow S/\left\langle g_{l}\right\rangle $ is $\left.  \omega
\right\vert _{\Delta_{1}}:\Delta_{1}\rightarrow\left\langle g_{l}\right\rangle
$ and $\Omega_{g_{l}}:H_{1}O((S/\left\langle g_{l}\right\rangle )\rightarrow
\mathbb{Z}_{2m}$ is given by the condition $\Omega_{g_{l}}(Y)=w$ iff
$\omega(y)=g_{l}^{w}$ where $\alpha(y)=Y$ and $\alpha:\Delta_{1}\rightarrow
H_{1}O((S/\left\langle g_{l}\right\rangle )$ is the abelianization.

\bigskip

As $g_{i}^{2}=f$, the number of fixed points of powers, rotation numbers
(except sign)\ and branched indices are the same for $g_{1}$ and $g_{2}$. Then
we have that the genus of $S/\left\langle g_{1}\right\rangle $ is equal to the
genus of $S/\left\langle g_{2}\right\rangle $, the signatures of $\Delta_{1}$
and $\Delta_{2}$ are also equal
\[
(k;-;[c_{1},...,c_{s}];\{-\})
\]
and
\[
\{\Omega_{g_{1}}(X_{j}(1))\}=\{\varepsilon_{j}\Omega_{g_{2}}(X_{j}(2))\}
\]
where $\{A_{i}(l),B_{i}(l),D_{q}(l),X_{j}(l)\}$ is the set of generators of
$H_{1}O((S/\left\langle g_{l}\right\rangle )$ given by canonical presentations
of $\pi_{1}O((S/\left\langle g_{l}\right\rangle )$, $l=1,2$. So in this case
the condition 1 in \ref{Classification} is verified automatically in order to
have the topological equivalence between $g_{1}$ and $g_{2}$.

\bigskip

The group $\Delta$ has a presentation of the form:

$\left\langle d_{1},d_{2},a_{i},b_{i},i=1,...,\frac{h-2}{2},x_{j}%
,j=1,...,r:d_{1}^{2}d_{2}^{2}\prod[a_{i},b_{i}]\prod x_{j}=1;x_{j}^{m_{j}%
}=1\right\rangle $ if $h$ is even and

$\left\langle d,a_{i},b_{i},i=1,...,\frac{h-1}{2},x_{j},j=1,...,r:d^{2}%
\prod[a_{i},b_{i}]\prod x_{j}=1;x_{j}^{m_{j}}=1\right\rangle $ if $h$ is odd.

\bigskip

\begin{lemma}
\label{Lemma 2}

If $h$ is even $\theta\circ\omega(d_{i})\notin\left\langle ab\right\rangle $,
and
\[
\left\langle \theta\circ\omega(a_{i}),\theta\circ\omega(b_{i}),\theta
\circ\omega(x_{j}),\theta\circ\omega(d_{1}d_{2})\right\rangle =\left\langle
ab\right\rangle .
\]

If $h$ is odd $\theta\circ\omega(d)\notin\left\langle ab\right\rangle $ and
\[
\left\langle \theta\circ\omega(a_{i}),\theta\circ\omega(b_{i}),\theta
\circ\omega(x_{j})\right\rangle =\left\langle ab\right\rangle .
\]

\end{lemma}

\textbf{Proof.}

Since $g_{1}$ and $g_{2}$ are anticonformal, the automorphism $g_{1}g_{2}$ is
conformal and $S/\left\langle g_{i}^{2}=f,g_{1}g_{2}\right\rangle $ is
orientable and uniformized by $\omega^{-1}(\left\langle f,g_{1}g_{2}%
\right\rangle )$.

Since $\omega^{-1}(\left\langle f,g_{1}g_{2}\right\rangle )$ has sign $+$ in
the signature, either $\omega(d_{i})$, where $i=1,2$, if $h$ is even, or
$\omega(d)$, if $h$ is odd, $\notin\omega^{-1}(\left\langle f,g_{1}%
g_{2}\right\rangle )=\omega^{-1}(\theta^{-1}\left\langle ab\right\rangle )$.
Hence $\theta\circ\omega(d_{1})$, $\theta\circ\omega(d_{2})$, $\theta
\circ\omega(d)$ are not in $\left\langle ab\right\rangle $.

If some element $y$ of $\left\langle a_{i},b_{i},x_{i},d_{1}d_{2}\right\rangle
$ satisfies that $\theta\circ\omega(y)$ is conjugate to $a$ or $b$, then
$d_{1}y$ or $dy$ belongs to $\omega^{-1}(\left\langle f,g_{1}g_{2}%
\right\rangle )$ that have all its elements orientation preserving, and this
is a contradiction, then $\theta\circ\omega(y)\in\left\langle ab\right\rangle
$.

If $h$ is odd and since $\theta\circ\omega(\Delta)=D_{n}$, this fact implies
\[
\left\langle \theta\circ\omega(a_{i}),\theta\circ\omega(b_{i}),\theta
\circ\omega(x_{j})\right\rangle =\left\langle ab\right\rangle
\]
In case $h$ even if $\theta\circ\omega(d_{1})=\theta\circ\omega(d_{2})$,
$\left\langle \theta\circ\omega(a_{i}),\theta\circ\omega(b_{i}),\theta
\circ\omega(x_{j})\right\rangle =\left\langle ab\right\rangle $ and if
$\theta\circ\omega(d_{1})\neq\theta\circ\omega(d_{2})$ then
\[
\left\langle \theta\circ\omega(d_{1}d_{2}),\theta\circ\omega(a_{i}%
),\theta\circ\omega(b_{i}),\theta\circ\omega(x_{j})\right\rangle =\left\langle
ab\right\rangle .
\]

$\square$

\bigskip

We shall consider separately two cases depending on the parity of the genus of
$S/\left\langle g_{1},g_{2}\right\rangle $.

\bigskip

\bigskip

\subsection{Case genus of $S/\left\langle g_{1},g_{2}\right\rangle =h$ even}

In this case, except when the genus of $S/\left\langle g_{i}\right\rangle $ is
two, the two square roots are topologically equivalent:

\begin{theorem}
\label{Thm1}Let $g_{i}$, $i=1,2$, be two anticonformal square roots of a
conformal automorphism $f$ of even order $m$. Assume that the genus of
$S/\left\langle g_{1},g_{2}\right\rangle $ is even and $\left\langle
g_{1},g_{2}\right\rangle $ has order $n\times2m$ with $n$ even. If the genus
of $S/\left\langle g_{i}\right\rangle $ is not two then $g_{1}$ and $g_{2}$
are topologically equivalent.
\end{theorem}

\bigskip

\textbf{Proof.}

We follow the notation fixed at the beginning of this Section, i.e. we have
$\mathbb{H}^{2}/\Delta\cong S/\left\langle g_{1},g_{2}\right\rangle $ and,
since the genus of $S/\left\langle g_{1},g_{2}\right\rangle $ is even,
$\Delta$ has presentation:%
\begin{equation}
\left\langle d_{1},d_{2},a_{i},b_{i},i=1,...,\frac{h-2}{2},x_{j}%
,j=1,...,r:d_{1}^{2}d_{2}^{2}\prod[a_{i},b_{i}]\prod x_{j}=1;x_{j}^{m_{j}%
}=1\right\rangle \tag{*}%
\end{equation}

As before $\omega:\Delta\rightarrow G=\left\langle g_{1},g_{2}\right\rangle $
is the monodromy of $S\rightarrow S/\left\langle g_{1},g_{2}\right\rangle $
and $\theta:G\rightarrow\left\langle g_{1},g_{2}\right\rangle /\left\langle
f\right\rangle =D_{n}$.

\bigskip

Since $\theta\circ\omega(\Delta)=D_{n}$ then $\theta\circ\omega\left\langle
a_{i},b_{i},x_{j},d_{1}d_{2}\right\rangle =\left\langle ab\right\rangle $, let
$y\in\left\langle a_{i},b_{i},x_{j}\right\rangle $ such that $\theta
\circ\omega(y)=ab$.

We call $\eta_{a}:D_{n}=\left\langle a,b\right\rangle \rightarrow\left\langle
a,b\right\rangle /\left\langle a\right\rangle =\mathcal{P}\{[0],...,[n-1]\}$,
where $[i]=(ab)^{i}\{1,a\}$, and given by the action of each element of
$D_{n}$ on the cosets $\left\langle a,b\right\rangle /\left\langle
a\right\rangle $:%
\[
\eta_{a}(z)[i]=[j]\text{, if }z(ab)^{i}\{1,a\}=(ab)^{j}\{1,a\}
\]
where $z\in D_{n}$. We have $\eta_{a}\circ\theta\circ\omega(y)[i]=[i+1]$. For
$\eta_{b}:\left\langle a,b\right\rangle \rightarrow\left\langle
a,b\right\rangle /\left\langle b\right\rangle $, analogously defined, we have
also $\eta_{b}\circ\theta\circ\omega(y)[i]=[i+1]$.

We shall contruct a set of generators of $Stab_{[0]}(\eta_{a}\circ\theta
\circ\omega)=\omega^{-1}(\left\langle g_{1}\right\rangle )=\pi_{1}%
O(S/\left\langle g_{1}\right\rangle )=\Delta_{1}$ and $Stab_{[0]}(\eta
_{b}\circ\theta\circ\omega)=\omega^{-1}(\left\langle g_{2}\right\rangle
)=\pi_{1}O(S/\left\langle g_{2}\right\rangle )=\Delta_{2}$. We apply the
Reidemeister-Schreier method (see for instance \cite{F}, Section 8) to obtain
a presentation of $\Delta_{l}$ from the presentation (*) of $\Delta$ and the
monodromy $\eta_{a}\circ\theta\circ\omega$. The orientation reversing
generators of $\Delta_{1}$ and $\Delta_{2}$ are:%
\[
y^{i}d_{1}(y^{i^{\prime}})^{-1},y^{i}d_{2}(y^{i^{\prime}})^{-1},\text{ where
}[i^{\prime}]=\eta_{a}(d_{i})[i]\text{, for }\Delta_{1}%
\]
and%
\[
y^{i}d_{1}(y^{i^{\prime\prime}})^{-1},y^{i}d_{2}(y^{i^{\prime\prime}}%
)^{-1},\text{ where }[i^{\prime\prime}]=\eta_{b}(d_{i})[i]\text{, for }%
\Delta_{2}%
\]

Let $\alpha_{l}:\Delta_{l}\rightarrow H_{1}O(S/\left\langle g_{l}\right\rangle
),l=1,2,$ be the abelianization morphisms and we will denote $D_{q}%
^{i}(1)=\alpha_{1}(y^{i}d_{q}(y^{i^{\prime}})^{-1})$ and $D_{q}^{i}%
(2)=\alpha_{2}(y^{i}d_{q}(y^{i^{\prime\prime}})^{-1})=D_{q}^{i}(2)$, $q=1,2$,
$i=0,...,n-1$.

Now we denote $N(1)$
\[
N(1)=\alpha_{1}(%
{\textstyle\prod\nolimits_{q}}
{\textstyle\prod\nolimits_{i}}
y^{i}d_{q}(y^{i^{\prime}})^{-1})=%
{\textstyle\sum\nolimits_{q=1,2}^{i=0,...,n-1}}
D_{q}^{i}(1).
\]
Analogously for $N(2)$.

In the same way we construct the orientation preserving generators for the
presentation of $\pi_{1}O(S/\left\langle g_{l}\right\rangle )$. If
$w\in\{a_{i},b_{i},x_{j}\}$ we have the generators $y^{i}w(y^{i^{\prime}%
})^{-1}$, where $\eta_{a}(w)[i]=i^{\prime}$. We denote $W^{i}(1)=\alpha
_{l}(y^{i}w(y^{i^{\prime}})^{-1})\,,$ $i=0,...,n-1$. Analogously
$W^{i}(2)=\alpha_{l}(y^{i}w(y^{i^{\prime\prime}})^{-1})$.

The elliptic elements of every canonical presentation of $\Delta_{l}=\pi
_{1}O(S/\left\langle g_{l}\right\rangle )$ are conjugate to powers of $x_{j}$.
Note that $\eta_{a}\circ\theta\circ\omega(x_{j})=\eta_{a}(ab)^{m/n_{j}}$ then
the order of $\eta_{a}\circ\theta\circ\omega(x_{j})$ is $n_{j}$ a divisor of
$n$, the elliptic generators of canonical presentations of $\pi_{1}%
O(S/\left\langle g_{l}\right\rangle )$ have the form $y^{i}x_{j}^{n_{j}}%
(y^{i})^{-1}$. The elliptic element $y^{i}x_{j}^{n_{j}}(y^{i})^{-1}$
corresponds to a period $c_{p}=m/n_{j}$, $p=1,...,s$, in the signature
$(k;-;[c_{1},...,c_{s}];\{-\})$ of $\Delta_{l}$.

Lifting the long relation of $\Delta$ we obtain $n$ relations for
$H_{1}O(S/\left\langle g_{l}\right\rangle )$ in the new generators
$\{D_{q}^{i}(l),W^{i}(l)\}$. Adding these $n$ relations we have a new relation
containing two times all the generators of the form $D_{q}^{i}(l)$, one time
each generator of the form $X_{j}^{i}(l),A_{j}^{i}(l)$ and $B_{j}^{i}(l)$ and
one time the opposite of each generator of the form $A_{j}^{i}(l)$ and
$B_{j}^{i}(l)$. Using the commutativity of $H_{1}O(S/\left\langle
g_{l}\right\rangle )$ we obtain:%
\[
2%
{\textstyle\sum\nolimits_{q=1,2}^{i=0,...,n-1}}
D_{q}^{i}(l)+%
{\textstyle\sum\nolimits_{j=1,...,r}^{i=0,...,n-1}}
X_{j}^{i}(l)=0
\]

Note that if $\theta\circ\omega(x_{j})$ have order $n_{j}$, where $n_{j}$
divides $m_{j}$, and $\prod\limits_{t=s_{1}}^{s_{n_{j}}}y^{t}x_{j}^{n_{j}%
}(y^{t})^{-1}$ is an elliptic generator of a canonical presentation of
$\Delta_{l}$. Then $%
{\textstyle\sum\nolimits_{j}^{i=0,...,n-1}}
X_{j}^{i}(l)$ is the sum of the homology classes determined by all conic
points of the orbifold $S/\left\langle g_{l}\right\rangle $ (all elliptic
generators of a canonical presentation of $\Delta_{l}$).

By the Lemma \ref{Lemma1} we have that there is a canonical generator system
$\{A_{i}(l),B_{i}(l),D_{q}^{\prime}(l),%
{\textstyle\sum\nolimits_{t=1}^{n_{j}}}
X_{j}^{s_{t}}(l)\}$, where $%
{\textstyle\sum\nolimits_{q=1,2}}
D_{q}^{\prime}(l)=%
{\textstyle\sum\nolimits_{q=1,2}^{i=0,...,n-1}}
D_{q}^{i}(l)=N(l)$.

Since $\omega(\prod\limits_{t=s_{1}}^{s_{n_{j}}}y^{t}x_{j}^{n_{j}}(y^{t}%
)^{-1})$ is an even power of $g_{i}$ and $g_{1}^{2}=g_{2}^{2}$,
\[
\omega(\prod\limits_{t=s_{1}}^{s_{n_{j}}}y^{t}x_{j}^{n_{j}}(y^{t})^{-1}%
)=g_{1}^{2w_{j}}=g_{2}^{2w_{j}}%
\]
then $\Omega_{g_{1}}(%
{\textstyle\sum\nolimits_{t=1}^{n_{j}}}
X_{j}^{s_{t}}(1))=2w_{j}=\Omega_{g_{2}}(%
{\textstyle\sum\nolimits_{t=1}^{n_{j}}}
X_{j}^{s_{t}}(2))$. If there are elements $X_{j}$ of order two and since
$S/\left\langle g_{i}\right\rangle $ has not genus $2$, then by the
classification in subsection \ref{Classification} the automorphisms $g_{1}$
and $g_{2}$ are topologically equivalent.

Assume now that there are no $X_{j}$ of order two, then to now if $g_{1}$ and
$g_{2}$ are topologically equivalent we need:%
\[
\Omega_{g_{1}}(%
{\textstyle\sum\nolimits_{q=1,2}}
D_{q}^{\prime}(1))=\Omega_{g_{2}}(%
{\textstyle\sum\nolimits_{q=1,2}}
D_{q}^{\prime}(2))
\]

Since $\left\langle g_{1}\right\rangle $ is abelian and $\omega(y^{i}%
d_{q}(y^{i^{\prime}})^{-1})\in\left\langle g_{1}\right\rangle $, we have that
\[
\omega(%
{\textstyle\prod\nolimits_{q}}
{\textstyle\prod\nolimits_{i}}
y^{i}d_{q}(y^{i^{\prime}})^{-1})=\omega(%
{\textstyle\prod\nolimits_{i}}
y^{i}d_{1}d_{2}(y^{i^{\prime}})^{-1})=\omega((d_{1}d_{2})^{n}).
\]
and it is an even power of $g_{l}$. Since $g_{1}^{2}=g_{2}^{2}$, $g_{l}%
^{2w}=\omega(%
{\textstyle\prod\nolimits_{i}}
y^{i}d_{1}d_{2}(y^{i})^{-1})=\omega((d_{1}d_{2})^{n})$, note that
$\omega(d_{1}d_{2})\in g_{i}^{2}$.

We have $\Omega_{g_{l}}(%
{\textstyle\sum\nolimits_{q=1,2}^{i=0,1}}
D_{q}^{i}(l))=2w$ then :%
\begin{gather*}
\Omega_{g_{1}}(%
{\textstyle\sum\nolimits_{q=1,2}}
D_{q}^{\prime}(1))=\Omega_{g_{1}}(%
{\textstyle\sum\nolimits_{q=1,2}^{i=0,1}}
D_{q}^{i}(1))=2w\\
=\Omega_{g_{2}}(%
{\textstyle\sum\nolimits_{q=1,2}^{i=0,1}}
D_{q}^{i}(2))=\Omega_{g_{2}}(%
{\textstyle\sum\nolimits_{q=1,2}}
D_{q}^{\prime}(2))
\end{gather*}

$\square$

\bigskip

The following example shows that the condition that the genus of
$S/\left\langle g_{l}\right\rangle $ is not two is necessary.

\bigskip

\begin{example}
\label{Ex1}
\end{example}

Consider $G=C_{16}\times C_{2}=\left\langle u:u^{16}=1\right\rangle
\oplus\left\langle t:t^{2}=2\right\rangle $ and $g_{1}=(u,1)$, $g_{2}=(u,t)$.
Let $\Delta$ be an NEC\ group with signature $(2;-;[2])$ and $\omega
:\Delta\rightarrow G$ given by $\omega(d_{1})=(u,1)$, $\omega(d_{2}%
)=(u^{3},t)$, $\omega(x)=(u^{8},1)$. Note that genus of $S/\left\langle
g_{l}\right\rangle $ is two and the signature of $\Delta_{l}$ is
$(2;-;[2,2])$, $l=1,2$.

A set of generators of a canonical presentation of $\Delta_{1}$ is
\[
\{d_{2}d_{1}d_{2}^{-1},d_{1}d_{2}^{2},d_{2}^{-2}d_{1}^{-1}xd_{1}d_{2}%
^{2},d_{2}xd_{2}^{-1}\}.
\]
Since $\omega(d_{2}d_{1}d_{2}^{-1})=g_{1}$ and if we call $\alpha(d_{2}%
d_{1}d_{2}^{-1})=D_{1}(1)$ and $\alpha(d_{1}d_{2}^{2})=D_{2}(1)$, we have
$\Omega_{g_{1}}(D_{1}(1))=1\operatorname{mod}16$.

Considering the canonical presentation of $\Delta_{2}$:%
\[
\{d_{1}^{2}d_{2}d_{1}^{-2},d_{1}d_{2}d_{1},d_{1}^{-1}d_{2}^{-1}d_{1}%
^{-2}xd_{1}^{2}d_{2}d_{1},x\},
\]
we have $\Omega_{g_{2}}(D_{1}(2))=3\operatorname{mod}16$.

Now $\Omega_{g_{1}}(D_{1}(1)+D_{2}(1))=8\operatorname{mod}16$ and
\[
\omega(d_{2}^{-2}d_{1}^{-1}xd_{1}d_{2}^{2})=g_{1}^{8},\omega(d_{2}xd_{2}%
^{-1})=g_{1}^{8},
\]
let $z=z_{1}=\gcd\{8,8,8,16\}=8$. We have $\Omega_{g_{1}}(D_{1}^{0}%
(1))\operatorname{mod}8=1\operatorname{mod}8$.

In the same way $\Omega_{g_{1}}(D_{1}^{0}(2))\operatorname{mod}%
8=3\operatorname{mod}8\neq\pm\Omega_{g_{1}}(D_{1}^{0}(1))\operatorname{mod}%
8=1\operatorname{mod}8$. Hence the automorphisms $g_{1}$ and $g_{2}$ are not
topologically equivalent. \ \ $\square$

\bigskip

\bigskip

\subsection{Case genus of $S/\left\langle g_{1},g_{2}\right\rangle =h$ odd}

To have an analogous result to Theorem \ref{Thm1} we need now the extra
condition that $\left\langle g_{1},g_{2}\right\rangle $ must be abelian:

\bigskip

\begin{theorem}
\label{Thm2}Let $g_{i}$, $i=1,2$, be two anticonformal square roots of a
conformal automorphism $f$ of even order $m$. Assume that genus of
$S/\left\langle g_{1},g_{2}\right\rangle $ is odd and $\left\langle
g_{1},g_{2}\right\rangle $ has order $n\times2m$ with $n$ even. If
$\left\langle g_{1},g_{2}\right\rangle $ is abelian then $g_{1}$ and $g_{2}$
are topologically equivalent.
\end{theorem}

\bigskip

\textbf{Proof.}

By the condition of $\left\langle g_{1},g_{2}\right\rangle $ abelian and
$\left\langle g_{1},g_{2}\right\rangle \neq\left\langle g_{i}\right\rangle $,
we have that $\left\langle g_{1},g_{2}\right\rangle /\left\langle g_{i}%
^{2}\right\rangle $ is $D_{2}$ and the order of $\left\langle g_{1}%
,g_{2}\right\rangle $ is $4m$.

Let $\Delta$ be an NEC\ group with signature:
\[
(1+2h;-;[m_{1},...,m_{r}];\{-\}).
\]
such that $S/\left\langle g_{1},g_{2}\right\rangle =\mathbb{H}^{2}/\Delta$.

We have a canonical presentation of $\Delta$ as follows:%
\begin{equation}
\left\langle d,a_{i},b_{i},x_{j}:d^{2}%
{\textstyle\prod}
[a_{i},b_{i}]%
{\textstyle\prod}
x_{j}=1;x_{j}^{m_{j}}=1\right\rangle \tag{**}%
\end{equation}

By Lemma \ref{Lemma 2}\ $\theta\circ\omega(d)\notin\left\langle
ab\right\rangle $\ and there is $y\in\left\langle a_{i},b_{i},x_{j}%
\right\rangle $ such that $\theta\circ\omega(y)=ab$. Without loss of
generality we suppose $\theta\circ\omega(d)=a$. Let $\eta_{a}$ be the natural
projection $\eta_{a}:D_{2}=\left\langle a,b:a^{2}=b^{2}=(ab)^{2}%
=1\right\rangle \rightarrow D_{2}/\left\langle a\right\rangle $ and analogously
$\eta_{b}:D_{2}\rightarrow D_{2}/\left\langle b\right\rangle $.

By Reidemeister-Schreier method the group $\omega^{-1}(\theta^{-1}(a))=\pi
_{1}O(S/\left\langle g_{1}\right\rangle )$ has a presentation with generators
$W^{i}(1)=\alpha(y^{i}w(y^{i^{\prime}})^{-1})$, where $w\in\{d,a_{i}%
,b_{i},x_{j}\}$ and $i^{\prime}=\eta_{a}\circ\theta\circ\omega(w)(i)$. Let
$\alpha_{1}:\pi_{1}O(S/\left\langle g_{1}\right\rangle )\rightarrow
H_{1}O(S/\left\langle g_{1}\right\rangle )$ be the abelianization and
$W^{i}(1)=\alpha(y^{i}w(y^{i^{\prime}})^{-1})$.

The orientation reversing generators produce in homology the generators
$D^{i}(1)$, and we call $N=%
{\textstyle\sum}
D^{i}(1)$. Note that the orientation reversing generators in $\pi
_{1}O(S/\left\langle g_{1}\right\rangle )$ are $y^{i}d(y^{i^{\prime}})^{-1}$,
$i=1,2$.

The elliptic generators of the canonical presentations of $\omega^{-1}%
(\theta^{-1}(a))=\pi_{1}O(S/\left\langle g_{1}\right\rangle )$ expressed in
the generators given by Reidemeister-Schreier are of the form
\[
y^{i}x_{j}(y^{i^{\prime}})^{-1}\text{ or }x_{j}^{2}%
\]
where $x_{j}$ is an elliptic generator of $\Delta$. Note that if $\theta
\circ\omega(x_{j})=ab$, $X_{j}$ gives only an elliptic generator of
$\omega^{-1}(\theta^{-1}(a))$, in homology $X_{j}^{0}+X_{j}^{1}=X_{j}%
^{\nu_{01}}$, if $\theta\circ\omega(x_{j})=1$, $X_{j}$ give two of such
generators, producing in homology $X_{j}^{0}=X_{j}^{\nu_{0}}$ and $X_{j}%
^{0}=X_{j}^{\nu_{1}}$. By the lifting of the long relation of $\Delta$ we
have
\[
2N+%
{\textstyle\sum}
X_{j}^{\nu}(1)=0.
\]

In the case where there are conic points of order two in $S/\left\langle
g_{i}\right\rangle $ the automorphisms $g_{1}$ and $g_{2}$ are equivalent by
the conditions in \ref{Classification}, in the contrary case, applying the
Lemma \ref{Lemma1}, we obtain a generator system
\[
\{A_{i}(1),B_{i}(1),D_{q}^{\prime}(1),X_{j}^{\nu}(1)\},
\]
that comes from a canonical presentation of $\pi_{1}O(S/\left\langle
g_{1}\right\rangle )$ and with
\[%
{\textstyle\sum\nolimits_{q=1,2}}
D_{q}^{\prime}(1)=N.
\]
Applying the commutativity of $\left\langle g_{1},g_{2}\right\rangle $ we
have $\omega(%
{\textstyle\prod}
y^{i}d(y^{i^{\prime}})^{-1})=\omega(d^{2})=g_{1}^{2}$, then $\Omega_{g_{1}}(%
{\textstyle\sum\nolimits_{q=1,2}}
D_{q}^{\prime}(1))=2$.

Analogously for $\pi_{1}O(S/\left\langle g_{2}\right\rangle )=\omega^{-1}%
(\theta^{-1}(b))$. Using that $n$ is even and $g_{1}^{2}=g_{2}^{2}$, we have
again $\Omega_{g_{2}}(%
{\textstyle\sum\nolimits_{q=1,2}}
D_{q}^{\prime}(2))=2$. Hence if genus of $S/\left\langle g_{i}\right\rangle $
$>2$ we have that the topological type of $g_{1}$ and $g_{2}$ is the same.

Let consider now the case where genus of $S/\left\langle g_{i}\right\rangle $
$=2$. We have obtained that $2=\Omega_{g_{i}}(%
{\textstyle\sum\nolimits_{q=1,2}}
D_{q}^{\prime}(i))$, then $\gcd\{\Omega_{g_{i}}(%
{\textstyle\sum\nolimits_{q=1,2}}
D_{q}^{\prime}(i)),\Omega_{g_{i}}(X_{v}(i))\}=2$. Since $\Omega_{g_{1}}%
(D_{1}^{\prime}(l))$ is odd, then $\Omega_{g_{1}}(D_{1}^{\prime}%
(1))=\Omega_{g_{2}}(D_{1}^{\prime}(2))\operatorname{mod}2$ and $g_{1}$ and
$g_{2}$ are topologically equivalent.

$\square$

\bigskip

The following example shows that the condition $\left\langle g_{1}%
,g_{2}\right\rangle $ abelian is necessary in the Theorem \ref{Thm2}:

\bigskip

\begin{example}
\label{Ex2}
\end{example}

Consider $G=\left\langle u,c:u^{2m}=c^{2}=1;cuc=u^{m+1}\right\rangle
=C_{2m}\rtimes C_{2}$, $m$ even $>2$. Note that the order of $G$ is $4m$.

Let $\Delta$ be an NEC\ group with signature $(1+2h,-,[m,2,2])$ and let
\[
\left\langle d,a_{i},b_{i},x_{j}:d^{2}%
{\textstyle\prod}
[a_{i},b_{i}]x_{1}x_{2}x_{3}=1;x_{1}^{m}=x_{2}^{2}=x_{3}^{2}=1\right\rangle
\]
be a canonical presentation of $\Delta$.

We construct the monodromy:%
\begin{align*}
\omega &  :\Delta\rightarrow G\\
d  &  \longmapsto u,x_{1}\longmapsto u^{-2},x_{2}\longmapsto c,x_{3}%
\longmapsto c,
\end{align*}
and $a_{i}$, $b_{i}$ can be sent to any element of $\left\langle
z^{2}\right\rangle <G$.

We call $g_{1}=u^{1+m/2}$ and $g_{2}=cu$. Then $g_{1}^{2}=u^{m+2}=g_{2}^{2}$.

1. Topological invariant of $g_{1}$ for condition 2 in \ref{Classification}.

The monodromy of the covering $\mathbb{H}^{2}/\omega^{-1}(\left\langle
g_{1}\right\rangle )\overset{2:1}{\rightarrow}\mathbb{H}^{2}/\Delta$ is given
by $\pi_{a}\circ\theta\circ\omega$, where $\pi_{a}:C_{2}\times C_{2}%
=\left\langle a,b\right\rangle \rightarrow C_{2}=\left\langle a,b\right\rangle
/\left\langle \theta\circ\omega(g_{1})=a\right\rangle =\{0,1;+\}$. We have:%
\[
\pi_{a}\circ\theta\circ\omega:d\longmapsto0,x_{1}\longmapsto0,x_{2}%
\longmapsto1,x_{3}\longmapsto1
\]

Then the signature of $\omega^{-1}(\left\langle g_{1}\right\rangle )$ is
$(t,-,[m,m])$ where $t$ is an even integer.

Orientation reversing generators for $\omega^{-1}(\left\langle g_{1}%
\right\rangle )$ given by Reidemeister-Schreier method:
\[
d,x_{2}dx_{2}^{-1},
\]
giving in homology:\ $\alpha(d)=D^{0}(1)$, $\alpha(x_{2}dx_{2}^{-1})=D^{1}(1)$.

We have $\Omega_{g_{1}}(D^{0}(1)+D^{1}(1))=2$, since $\omega(dx_{2}dx_{2}%
^{-1})=u^{m+2}=(g_{1})^{2}$

2. Topological invariant of $g_{2}$ for condition 2 in \ref{Classification}.

The monodromy of the covering $H^{2}/\omega^{-1}(\left\langle g_{2}%
\right\rangle )\overset{2:1}{\rightarrow}H^{2}/\Delta$ is given by%
\[
\pi_{b}\circ\theta\circ\omega:d\longmapsto1,x_{1}\longmapsto0,x_{2}%
\longmapsto1,x_{3}\longmapsto1
\]

As before the signature of $\omega^{-1}(\left\langle g_{2}\right\rangle )$ is
$(t,-,[m,m])$. The orientation reversing generators for $\omega^{-1}%
(\left\langle g_{2}\right\rangle )$ are now
\[
x_{2}d,dx_{2}^{-1}%
\]
Giving in homology:\ $D^{0}(2)$, $D^{1}(2)$. Hence $\Omega_{g_{2}}%
(D^{0}(2)+D^{1}(2))\neq2$, since $\omega(x_{2}ddx_{2}^{-1})=u^{2}\neq
(g_{2})^{2}$.

Then $g_{1}$ and $g_{2}$ are not topologically equivalent. $\square$

\bigskip

From the above theorems we have:

\begin{corollary}
Let $g_{i}$, $i=1,2$, be two anticonformal square roots of a conformal
automorphism $f$ of order $m$ even. If genus of $S/\left\langle g_{i}%
\right\rangle $ is not two and $\left\langle g_{1},g_{2}\right\rangle $ is
abelian then $\left\langle g_{1}\right\rangle $ and $\left\langle
g_{2}\right\rangle $ are topologically equivalent.
\end{corollary}

\bigskip

Antonio F. Costa

Departamento de Matem\'{a}ticas Fundamentales 

Facultad de Ciencias, UNED

C. Juan del Rosal, 10

28040 Madrid, Spain
\end{document}